\documentclass[10pt]{article}
\usepackage{latexsym}
\usepackage{amsfonts}
\usepackage{enumerate}
\usepackage{multicol}

\title{The ``Thirty-seven Percent Rule'' and the Secretary Problem with Relative Ranks \\[.4in]}

\author{B\'{e}la Bajnok \\[.1in] {\small Department of Mathematics, Gettysburg College} \\
{\small 300 N. Washington Street, Gettysburg, PA 17325-1486 USA} \\{\small E-mail:  bbajnok@gettysburg.edu} \\ [.2in]
and \\[.2in]
Svetoslav Semov \\[.1in] {\small Department of Economics, Boston University} \\
{\small 270 Bay State Rd, Boston, MA 02215 USA} \\ {\small E-mail: semov@bu.edu}   \\
[.4in]}

\date{March 12, 2013}

\newtheorem{thm}{Theorem}

\begin{document}

\maketitle

\begin{abstract}

We revisit the problem of selecting an item from $n$ choices that appear before us in random sequential order so as to minimize the expected rank of the item selected.  In particular, we examine the stopping rule where we reject the first $k$ items and then select the first subsequent item that ranks lower than the $l$-th lowest-ranked item among the first $k$.  We prove that the optimal rule has $k \sim n/{\mathrm e}$, as in the classical secretary problem where our sole objective is to select the item of lowest rank; however, with the optimally chosen $l$, here we can get the expected rank of the item selected to be less than any positive power of $n$ (as $n$ approaches infinity).  We also introduce a common generalization where our goal is to minimize the expected rank of the item selected, but this rank must be within the lowest $d$.


{\em 2010 AMS codes:} Primary: 62L99; Secondary: 60C05.

\end{abstract}

\section{Introduction}

Suppose that we have a job opening for which we have $n$ applicants of whom we must hire exactly one.  The applicants appear before us in random order one at a time; all $n!$ permutations are equally likely.  We assume that we do not know anything about the applicants until we interview them but that we have a clear preference between any two whom we have already met with.  After interviewing a candidate, we can either accept them, in which case the search is over, or we reject them and move on to the next candidate.  Our hiring decision is constrained by the rule that rejected candidates cannot be recalled; therefore, terminating the process early and waiting too long both come with risks. (In particular, if we reject each of the first $n-1$ applicants then we are forced to hire the last one.)  What shall we do?

The answer, of course, depends on what our goal is exactly.  In the classical ``secretary problem,''  we aim to maximize the probability that we hire the best candidate.  It is easy to see that the best strategy is to reject the first $k$ candidates for some $k$ and then hire the applicant whom we prefer over all the ones that we have seen thus far; it was first proved by Gilbert and Mosteller in 1966 \cite{GilMos:1966a} --- and now well known --- that this stopping rule ${\cal R}_n(k)$ is optimal at $k \sim n/{\mathrm e}$.  The general appeal of the secretary problem can be partly attributed to this surprisingly attractive answer, known popularly as the ``thirty-seven percent rule.''  There is a vast literature and many variations of the secretary problem; see, for example, Freeman \cite{Fre:1983a}, Ferguson \cite{Fer:1989a}, Pfeifer \cite{Pfe:1989a}, Bruss and Ferguson \cite{BruFer:1993a}, Assad and Samuel-Cahn \cite{AssSam:1996a}, Quine and Law \cite{QuiLaw:1996a}, Krieger and Samuel-Cahn \cite{KriSam:2009a}, and their references.   

Rather than focusing solely on the best applicant, one might instead aim to minimize the expected rank of the applicant hired.  As we can see from Theorem \ref{En(k,l)} below, the stopping rule ${\cal R}_n(k)$ just described yields a rather high expected rank, approximately $n/2\mathrm{e}$.   Instead, as  proved by Bearden \cite{Bea:2006a} in 2006, with the optimally chosen $k \sim \sqrt{n}$, the expected rank of the chosen candidate is considerably less, only about $\sqrt{n}$.  A markedly better result can be achieved if we allow an adaptive strategy where our stopping rule is not given in advance but is dependent on the relative ranks of the applicants we have already seen.  Such an optimal strategy was found by Lindley \cite{Lin:1961a} in 1961; a few years later Chow et al. \cite{Cho:1964a} proved that, as $n$ approaches infinity, the expected rank of the hired candidate using Lindley's stopping rule tends to the constant value $$\prod_{k=1}^{\infty} \left( 1+ \frac{2}{k} \right)^{1/(k+1)} \approx 3.87.$$  This is an amazing result, though the strategy that achieves it can only be stated implicitly and employed via dynamic programming.

In Section 2 of this paper, we discuss an explicit and {\em a priori} stopping rule that generalizes Bearden's result.  Namely, we investigate the stopping rule ${\cal R}_n(k,l)$ that rejects the first $k$ candidates and then hires the first one after that that ranks lower than the $l$-th best candidate among the first $k$.  As we shall see, for fixed $l$, the optimal value of $k$ is about $$k \sim  n^{\frac{l}{l+1}}.$$  (Note that for $l=1$ we get Bearden's result.)  As it turns out, with this $k$, the value of $l$ for which the stopping rule ${\cal R}_n(k,l)$  minimizes the expected rank of the applicant hired is about $\log n -1$.  This is a particularly pleasing result as $$n^{\frac{\log n -1}{\log n}}=n/\mathrm{e};$$  in other words, we are rejecting the first $\sim n/\mathrm{e}$ applicants, just like in the classical secretary problem, but now the expected value of the rank of the applicant hired is substantially lower, about $\mathrm{e} \log n /2$.

In Section 3 of our paper we introduce a common generalization of the classical secretary problem (where we are only interested in the best applicant) and the variation just discussed (where we have no absolute requirements on the applicant hired).  Namely, we investigate what happens when we still want to minimize the rank of the applicant hired, but we insist on hiring one of the best $d$ applicants ($1 \leq d \leq n$).  The case $d=1$ then yields the classical secretary problem, and the case $d=n$ corresponds to situation we analyzed in Section 2; we also provide a full analysis of the case $d=2$. 
 
{\bf Acknowledgments.}  We thank Art Benjamin, Darren Glass, Benjamin Kennedy, Chuck Wessell, and Rebecca Zabel for valuable discussions of various parts of this paper.

\section{The expected rank resulting from strategy ${\cal R}_n(k,l)$}

Let $n$, $k$, and $l$ be fixed positive integers with $l \leq k \leq n-1$.  By the stopping rule ${\cal R}_n(k,l)$ we mean the strategy that selects an element $s$ (for {\em selected value}) of a permutation $a_1,a_2,\dots,a_n$ of $\{1,2,\dots,n\}$ as follows.  We let $t$ (for {\em test value}) denote the $l$-th lowest value among $a_1,\dots,a_k$ and set 
$$I=\{i \in \mathbb{N} \; | \; i \leq n-k \; \mbox{and} \; a_{k+i} < t\}.$$ We then let 
$$j = \left\{
\begin{array}{cl}
\min I & \mbox{if $I \neq \emptyset,$}\\
n-k & \mbox{if $I = \emptyset;$}\\
\end{array}\right.$$ 
the rule then selects $s=a_{k+j}$.  In other words, the first $k$ candidates are rejected, after which the first value, if there is one, is selected that ranks lower than the $l$-th best candidate among the first $k$, and the last value is selected otherwise.  We can determine the exact expected value of $s$ as follows.

\begin{thm} \label{En(k,l)}

The expected value $E_n(k,l)$ of the rank of the candidate hired following strategy ${\cal R}_n(k,l)$ is $$E_n(k,l)=\frac{n+1}{2} \cdot \left( \frac{l}{k+1} + \frac{{n-l \choose k-l}}{{n \choose k}} \right).$$

\end{thm}

{\em Proof.}   First we introduce some terminology.  We say that the sequence $a_1, a_2, \dots, a_n$ is {\em successful} if at least one of $a_{k+1}, a_{k+2}, \dots, a_n$ is less than $t$ (that is, the set $I$ above is nonempty); otherwise, if each of $a_{k+1}, a_{k+2}, \dots, a_n$ is more than $t$, we say that the sequence is {\em unsuccessful}.  We will consider the cases of successful and unsuccessful sequences separately.

Suppose first that our search is successful and we select the value $s=a_{k+j}$ for some $j \in \mathbb{N}$.  Note that 
$$1 \leq s \leq t-1;$$ furthermore, since exactly $l$ of the values $a_1,a_2,\dots,a_{k+j}$ are less than $t$ and $k-l+j-1$ are more than $t$, we must have
$l+1 \leq t \leq n-k+l-j+1$.  Consequently, 
$$1 \leq j \leq n-k+l-t+1$$ and
$$l+1 \leq t \leq n-k+l$$ must hold. 

We need to determine, for fixed values of $t$, $j$, and $s$,  the number of sequences $a_1, a_2, \dots, a_n$ for which
\begin{itemize}

\item  $t$ is the $l$-th lowest value among $a_1,a_2,\dots,a_k$;
\item  $a_{k+1},a_{k+2},\dots,a_{k+j-1}$ are each more than $t$; and
\item  $s=a_{k+j}$ is less than $t$. 

\end{itemize} We do this as follows.

First, we choose the sequence $a_1,a_2,\dots,a_k$.  Since $l-1$ of the $k$ elements are less than $t$ but unequal to $s$ and $k-l$ of them are more than $t$, we have $${t-2 \choose l-1} \cdot {n-t \choose k-l}$$ ways to select the values; we then have $k!$ ways to arrange them into a sequence.

Next, we choose the sequence $a_{k+1}, a_{k+2}, \dots, a_{k+j-1}$.  Since each of these terms must be more than $t$ and distinct from those $k-l$ values among $a_1,a_2,\dots,a_k$ that are also more than $t$, we have 
$${n-k+l-t \choose j-1}$$ ways to select the values; we then have $(j-1)!$ ways to arrange them into a sequence.     

Finally, we have $(n-k-j)!$ ways to arrange the remaining $n-k-j$ elements of $\{1,2,\dots,n\}$ into the sequence $a_{k+j+1}, a_{k+j+2}, \dots, a_n$.

In summary, there are
$${t-2 \choose l-1} \cdot {n-t \choose k-l} \cdot {n-k+l-t \choose j-1} \cdot k! \cdot (j-1)! \cdot (n-k-j)!$$
sequences $a_1, a_2, \dots, a_n$ satisfying the three conditions above, which we may rewrite as
$${t-1 \choose l-1} \cdot {n-t \choose k-l} \cdot {n-k-j \choose t-l-1} \cdot \frac{(t-l) \cdot k! \cdot (n-k+l-t)!}{t-1}.$$  (Note that, in the case of a successful sequence, $t \neq 1$.)  Therefore, for a given pair of integers $t$ and $s$, with $1 \leq s \leq t-1$, the probability $P^{\surd}_n(k,l;t,s)$ that the stopping rule ${\cal R}_n(k,l)$ results in the (successful) selection of $s$ is
$$P^{\surd}_n(k,l;t,s)= \sum_{j=1}^{n-k+l-t+1} {t-1 \choose l-1} \cdot {n-t \choose k-l} \cdot {n-k-j \choose t-l-1} \cdot \frac{(t-l) \cdot k! \cdot (n-k+l-t)!}{(t-1) \cdot n!}.$$

We can simplify this expression by observing that ${n-k-j \choose t-l-1}$ enumerates the $(t-l)$-subsets of the set $\{1,2,\dots,n-k\}$ whose smallest element is $j$.  Letting $j$ vary from 1 to $n-k+l-t+1$ includes all possibilities (when $j=n-k+l-t+1$, our set consists of the greatest $t-l$ elements of the set $\{1,2,\dots,n-k\}$).  Therefore,
$$\sum_{j=1}^{n-k+l-t+1} {n-k-j \choose t-l-1} = {n-k \choose t-l}=\frac{(n-k)!}{(n-k+l-t)! \cdot (t-l)!},$$ yielding
\begin{eqnarray} \label{Psurd}
P^{\surd}_n(k,l;t,s)= \frac{1}{{n \choose k}} \cdot {t-1 \choose l-1} \cdot {n-t \choose k-l} \cdot \frac{1}{t-1}.
\end{eqnarray}

Let us turn now to the case of unsuccessful sequences.  The sequence $a_1, a_2, \dots, a_n$ leads to an unsuccessful search when each of $a_{k+1}, a_{k+2}, \dots, a_n$ is more than the $l$-th smallest element among $a_1, a_2, \dots, a_k$.  In other words, a sequence is unsuccessful exactly when each of $1,2,\dots,l$ is among the first $k$ terms of the sequence.  Therefore, for a given integer $s$, with $l+1 \leq s \leq n$, the probability $P^{X}_n(k,l;t,s)$ that the stopping rule ${\cal R}_n(k,l)$ results in the (unsuccessful) selection of $s$ is    
\begin{eqnarray} \label{PX}
P^{X}_n(k,l;s)= \frac{k!}{(k-l)!} \cdot \frac{(n-l-1)!}{n!}=\frac{1}{{n \choose k}} \cdot {n-l \choose k-l} \cdot \frac{1}{n-l}.
\end{eqnarray}

We can then exhibit the expected value of $s$ as
$$E_n(k,l) = \sum_{t=l+1}^{n-k+l} \sum_{s=1}^{t-1} P^{\surd}_n(k,l;t,s) \cdot s + 
    \sum_{s=l+1}^{n} P^{X}_n(k,l;s) \cdot s,$$
where
\begin{eqnarray*}
\sum_{t=l+1}^{n-k+l} \sum_{s=1}^{t-1} P^{\surd}_n(k,l;t,s) \cdot s & = &   \sum_{t=l+1}^{n-k+l} \frac{1}{{n \choose k}} \cdot {t-1 \choose l-1} \cdot {n-t \choose k-l} \cdot \frac{1}{t-1} \cdot \frac{t^2-t}{2} \\
& = &  \frac{l }{2 {n \choose k}} \cdot \sum_{t=l+1}^{n-k+l}  {n-t \choose k-l} \cdot {t \choose l}
\end{eqnarray*}
and
\begin{eqnarray*}
\sum_{s=l+1}^{n} P^{X}_n(k,l;s) \cdot s  & = &  \frac{1}{{n \choose k}} \cdot {n-l \choose k-l} \cdot \frac{1}{n-l} \cdot \left(\frac{n^2+n}{2}-\frac{l^2+l}{2} \right) \\
& = &  \frac{n+l+1 }{2 {n \choose k}} \cdot {n-l \choose k-l}.
\end{eqnarray*}

Next, observe that
$$\sum_{t=l}^{n-k+l}  {n-t \choose k-l} \cdot {t \choose l}={n+1 \choose k+1},$$ since both sides of the equation count the number of $(k+1)$-subsets of the set $\{1,2,\dots,n+1\}$.  Indeed, the number of $(k+1)$-subsets of $\{1,2,\dots,n+1\}$ whose $(l+1)$-st smallest element equals $t+1$ is $${ n-t \choose k-l}\cdot  {t \choose l};$$ as $t$ ranges from $l$ to $n-k+l$, we get all $(k+1)$-subsets of $\{1,2,\dots,n+1\}$.  Therefore,
$$\sum_{t=l+1}^{n-k+l}  {n-t \choose k-l} \cdot {t \choose l}={n+1 \choose k+1}-{n-l \choose k-l},$$ and thus
\begin{eqnarray*}
 E_n(k,l) & = & \frac{l }{2 {n \choose k}} \cdot \left( {n+1 \choose k+1}-{n-l \choose k-l}    \right) + \frac{n+l+1 }{2 {n \choose k}} \cdot {n-l \choose k-l} \\
& = & \frac{n+1}{2} \cdot \left( \frac{l}{k+1} + \frac{{n-l \choose k-l}}{{n \choose k}} \right), 
\end{eqnarray*}
as claimed.  $\Box$

We can now use Theorem \ref{En(k,l)} to find the stopping rule ${\cal R}_n(k,l)$ that minimizes the expected rank $E_n(k,l)$.  For $l=1$ we have
$$E_n(k,1)=\frac{n+1}{2} \cdot \left( \frac{1}{k+1} + \frac{k}{n} \right);$$ this quantity attains its minimum value when $k =\sqrt{n}-1$, and thus the stopping rule ${\cal R}_n(k,1)$ is optimal when $k=\lfloor \sqrt{n}-1 \rfloor$ or $k=\lceil \sqrt{n}-1 \rceil$, confirming Proposition 1 in \cite{Bea:2006a} (although the result there was under slightly different assumptions).

Minimizing for $l=2$ exactly is more complicated.  We have 
$$E_n(k,2)=\frac{n+1}{2} \cdot \left( \frac{2}{k+1} + \frac{k^2-k}{n^2-n} \right),$$ and this quantity attains its minimum value when 
$$k =\frac{1}{2} \cdot \left( \left(\sqrt{(2n-1)^4-1}+(2n-1)^2 \right)^{1/3} -1 + \left(\sqrt{(2n-1)^4-1}+(2n-1)^2 \right)^{-1/3} \right);$$ for large $n$ this can be approximated as $k \sim n^{2/3}$.  

More generally, for fixed $l$ and for large $k$ and $n$, we have 
$$E_n(k,l) = \frac{n+1}{2} \cdot \left( \frac{l}{k+1} + \frac{k \cdot (k-1) \cdot \cdots \cdot (k-l+1)}{n \cdot (n-1) \cdot \cdots \cdot (n-l+1)}\right) \sim \frac{n}{2} \cdot \left( \frac{l}{k} + \left( \frac{k}{n}\right)^l \right),$$ which attains its minimum at $k=n^{l/(l+1)}$.  For this value of $k$, we have
$$E_n(k,l) \sim \frac{n}{2} \cdot \left( \frac{l}{k} + \left( \frac{k}{n}\right)^l \right) = \frac{n}{2} \cdot \left( \frac{l}{n^{l/(l+1)}} + \left( \frac{n^{l/(l+1)}}{n}\right)^l \right) = \frac{l+1}{2} \cdot n^{1/(l+1)},$$ and this is minimal when $l=\log n -1$, in which case we have $$k=n^{l/(l+1)}= \frac{n}{\mathrm{e}} $$  and 
$$E_n(k,l) \sim \frac{l+1}{2} \cdot n^{1/(l+1)} = \frac{\mathrm{e}}{2}   \log n .$$

We can thus see that whether our goal is to maximize the probability of hiring the top candidate, as in the classical secretary problem, or to minimize the expected rank of the applicant hired, we should let the first $k \sim n/\mathrm{e}$ applicants pass; but by hiring the $l$-th best candidate after the first $k$ with $l \sim \log n -1$ yields an expected rank of approximately $\mathrm{e} \log n /2$, a substantial improvement compared to the classical rule that aims to select the best candidate ($l=1$).

\section{A generalization}

In the classical secretary problem, one's goal is to select the absolute lowest ranked applicant, while in the problem we just considered, we make no demands on the rank of the applicant selected.  As a common generalization, here we consider the case when we still want to minimize the expected rank of the applicant hired, but we are only interested when this rank is within the top $d$ (for a given $1 \leq d \leq n$).  We may reformulate the problem in terms of rewards: If an applicant of absolute rank $s$ is hired, our reward is $v_{n,d}(s)$ where     
$$v_{n,d}(s) = \left\{
\begin{array}{cl}
n+1-s & \mbox{if $s \leq d,$}\\
0 & \mbox{otherwise.}\\
\end{array}\right.$$ 
The classical secretary problem then corresponds to the case of $d=1$ (when the reward is $n$ if we hire the best applicant but 0 otherwise), and the no-demand problem of the previous section corresponds to $d=n$ (where the reward is $n$ for the best applicant, $n-1$ for the second best, and so on).

We can use our previous expressions (\ref{Psurd}) and (\ref{PX}) to compute the expected reward $V_{n,d}(k,l)$ following our stopping rule ${\cal R}_n(k,l)$ as
\begin{eqnarray} \label{Vformula}
V_{n,d}(k,l)=\sum_{t=l+1}^{n-k+l} \sum_{s=1}^{\min\{t-1,d\}} P^{\surd}_n(k,l;t,s) \cdot (n+1-s) + 
    \sum_{s=l+1}^{d} P^{X}_n(k,l;s) \cdot (n+1-s).
\end{eqnarray}

We can then easily see that $$V_{n,n}(k,l) = (n+1) - E_n(k,l).$$  Indeed, we have
\begin{eqnarray*} 
E_n(k,l)  + V_{n,n}(k,l) & = & \sum_{t=l+1}^{n-k+l} \sum_{s=1}^{t-1} P^{\surd}_n(k,l;t,s) \cdot s + 
    \sum_{s=l+1}^{n} P^{X}_n(k,l;s) \cdot s \\
& & + \sum_{t=l+1}^{n-k+l} \sum_{s=1}^{t-1} P^{\surd}_n(k,l;t,s) \cdot (n+1-s) + 
    \sum_{s=l+1}^{n} P^{X}_n(k,l;s) \cdot (n+1-s) \\ \\
& = & (n+1) \cdot \left( \sum_{t=l+1}^{n-k+l} \sum_{s=1}^{t-1} P^{\surd}_n(k,l;t,s) + 
    \sum_{s=l+1}^{n} P^{X}_n(k,l;s) \right) \\ \\
& = & n+1,
\end{eqnarray*}
since
$$\sum_{t=l+1}^{n-k+l} \sum_{s=1}^{t-1} P^{\surd}_n(k,l;t,s) + 
    \sum_{s=l+1}^{n} P^{X}_n(k,l;s)$$
is the sum of all probabilities and thus equals 1.

As we have seen in the previous section, with $k$ and $l$ approximately $n/\mathrm{e}$ and $\log n -1$, respectively, we have $$E_n(k,l) \sim \frac{\mathrm{e}}{2}   \log n,$$ and therefore  
$$\lim_{n \to \infty} \frac{\max\{V_{n,n}(k,l) \mid k, l \}}{n} = \lim_{n \to \infty} \frac{\min\{n+1-E_{n}(k,l) \mid k, l \}}{n} =1.$$  Thus we can rest assured: with a large pool of applicants, when following a simple strategy but not deeming any applicant unacceptable, we are poised to end up with one of the best secretaries anyway.

In general, we are interested in finding (or estimating) for every (fixed) $d$, the value of
$$c_d= \lim_{n \to \infty} \frac{\max\{V_{n,d}(k,l) \mid k, l \}}{n}.$$

Computational data suggest that
$$c_1 \approx 0.37, \; c_2 \approx 0.51, \; c_3 \approx 0.63, \; c_4 \approx 0.71,  \; c_5 \approx 0.77, \; c_6 \approx 0.81, \; c_7 \approx 0.84, \; c_8 \approx 0.87, \dots .$$  According to the classical secretary problem, and as we confirm below, we have $c_1=1/\mathrm{e};$  we will also find the exact value of $c_2$.

First, we establish the following explicit results.

\begin{thm} \label{Vnd}

Let $H_n(k)$ denote the harmonic number $\sum_{i=1}^{n-k} 1/(n-i)$.  The expected rewards $ V_{n,d}(k,l)$ for the stopping rule ${\cal R}_n(k,l)$ in the cases of $d=1$ and $d=2$ are as follows.
\begin{eqnarray}
V_{n,1}(k,l) & = & \left\{
\begin{array}{ll}
k \cdot (H_n-H_k) & \mbox{if $l=1$}\\ \\  
\frac{n}{l-1} \cdot \left( \frac{k}{n} -  \frac{{n-l \choose k-l}}{{n \choose k}} \right) & \mbox{if $l \geq 2$}\\
\end{array}\right. \\ \nonumber \\  \nonumber \\
V_{n,2}(k,l) & = & \left\{
\begin{array}{ll}
\frac{2n-1}{n} \cdot k \cdot (H_n-H_k) - \frac{k}{n} \cdot (n-k-1) & \mbox{if $l=1$}\\ \\
\frac{2n-1}{l-1} \cdot \left( \frac{k}{n} -  \frac{{n-l \choose k-l}}{{n \choose k}} \right) & \mbox{if $l \geq 2$}\\
\end{array}\right.
\end{eqnarray}

\end{thm}

{\em Proof.}  Assume first that $l \geq d$.  Then  (\ref{Vformula}) becomes 
$$V_{n,d}(k,l)=\sum_{t=l+1}^{n-k+l} \sum_{s=1}^{d} P^{\surd}_n(k,l;t,s) \cdot (n+1-s);$$
substituting (\ref{Psurd}) and simplifying yields
\begin{eqnarray} \label{Vl>=d}
V_{n,d}(k,l)=\frac{d \cdot (2n+1-d)}{2 {n \choose k}} \cdot  \sum_{t=l+1}^{n-k+l}{t-1 \choose l-1} \cdot {n-t \choose k-l} \cdot \frac{1}{t-1}.
\end{eqnarray}
With $l=d=1$, (\ref{Vl>=d}) becomes
\begin{eqnarray} \label{Vl=d=1}
V_{n,1}(k,1)=\frac{n}{{n \choose k}} \cdot  \sum_{t=2}^{n-k+1} {n-t \choose k-1} \cdot \frac{1}{t-1}.
\end{eqnarray}
Our result for $V_{n,1}(k,1)$ then follows from the identity
\begin{eqnarray} \label{Stirlingl=1}
\sum_{t=2}^{n-k+1} {n-t \choose k-1} \cdot \frac{1}{t-1} = 
{n-1 \choose k-1} \cdot (H_n-H_k),
\end{eqnarray}
which we prove next.

Note that (\ref{Stirlingl=1}) is equivalent to 
\begin{eqnarray} \label{Stirlingl=1'}
\sum_{t=2}^{n-k+1} {n-t \choose k-1} \cdot (k-1)! \cdot {n-k \choose n-t+1-k} \cdot (n-t+1-k)! \cdot (t-2)! +(n-1)! \cdot H_k=
(n-1)! \cdot H_n.
\end{eqnarray}
Here we recognize that the right-hand side of (\ref{Stirlingl=1'}) is Stirling's cycle number  $\left[{n \atop 2}\right]$, which counts the number of arrangements of the first $n$ positive integers into two disjoint nonempty cycles.  (The order of the two cycles does not matter; as customary, we assume that the one containing 1 appears first.  Furthermore, we assume that the cycles are both listed so that their smallest element appears first.)  Indeed, if the first cycle has length $n-t$, then we have $$\frac{(n-1)!}{t!}$$ choices for it, leaving $(t-1)!$ choices for the second cycle.  As $t$ ranges from 1 to $n-1$, we get all possibilities, and thus 
$$\left[{n \atop 2}\right]=\sum_{t=1}^{n-1} \frac{(n-1)!}{t!} \cdot (t-1)!=\sum_{t=1}^{n-1} \frac{(n-1)!}{t}=(n-1)! \cdot H_n.$$ 

We prove (\ref{Stirlingl=1'}) by showing that 
\begin{eqnarray} \label{Stirlingl=1'Left}
\sum_{t=2}^{n-k+1} {n-t \choose k-1} \cdot (k-1)! \cdot {n-k \choose n-t+1-k} \cdot (n-t+1-k)! \cdot (t-2)! 
\end{eqnarray}
counts those cycle-decompositions that contain each of $1,2,\dots,k$ in the same (first) cycle, while $$(n-1)! \cdot H_k$$ counts all others.    

To verify the first claim, imagine that the first cycle (the one containing $1,2,\dots,k$) has length $n-t+1$; letting $t$ vary from 2 to $n-k+1$ assures that we have considered all possibilities where neither cycle is empty.  To count the number of choices for the first cycle, recall that it starts with 1; we then have  
$${n-t \choose k-1}$$ ways to choose places for $2,3,\dots,k$, and $(k-1)!$ ways to arrange them.  We also have $${n-k \choose n-t+1-k}$$ ways to choose the remaining elements in the first cycle, with $(n-t+1-k)!$ ways to arrange them.  Our second cycle has length $t-1$, so we have $(t-2)!$ choices with the smallest element appearing first there too.  Thus, the number of arrangements of $1,2,\dots, n$ into two disjoint nonempty cycles with each of $1,2,\dots,k$ in the first cycle is given by (\ref{Stirlingl=1'Left}).

To enumerate the remaining cycle-decompositions, start by arranging $1,2,\dots,k$ into two disjoint nonempty cycles: there are      
$\left[{k \atop 2}\right]$ ways to do this.  Then, place the elements $k+1,k+2,\dots,n$ into one of the two cycles one at a time; each can be put to the right of an already-placed element.  (The two leading terms in the cycles must remain.)  The number of choices to place these elements is $$\frac{(n-1)!}{(k-1)!},$$ proving our second claim as
$$ \left[{k \atop 2}\right] \cdot \frac{(n-1)!}{(k-1)!} = (n-1)! \cdot H_k.$$ This completes the proof of the formula for $V_{n,1}(k,1)$.

Next, we prove that for $l \geq 2$ we have
\begin{eqnarray} \label{Stirlingl>=2}
\sum_{t=l+1}^{n-k+l}{t-1 \choose l-1} \cdot {n-t \choose k-l} \cdot \frac{1}{t-1} = 
\frac{1}{l-1} \cdot \left( {n-1 \choose k-1} - {n-l \choose k-l} \right) ,
\end{eqnarray}
which then, with (\ref{Vl>=d}), establishes our results for $V_{n,1}(k,l)$ and $V_{n,2}(k,l)$ when $l \geq 2$.

We start the proof of (\ref{Stirlingl>=2}) by rewriting it as
\begin{eqnarray*}
\sum_{t=l}^{n-k+l}{t-1 \choose l-1} \cdot {n-t \choose k-l} \cdot \frac{1}{t-1} = 
\frac{1}{l-1} \cdot {n-1 \choose k-1} ,
\end{eqnarray*} and then as
\begin{eqnarray} \label{Stirlingl>=2'}
\sum_{t=l}^{n-k+l}{t-2 \choose l-2} \cdot (l-2)! \cdot {n-k \choose t-l} \cdot (t-l)! \cdot {n-t \choose k-l}\cdot (k-l)! \cdot (n-t-k+l)! = 
\frac{(n-1)! \cdot (k-l)! \cdot (l-2)!}{(k-1)!}.
\end{eqnarray}

In a manner similar to the one above, we establish this identity by showing that both sides count the number of arrangements of $1,2,\dots,n$ into two disjoint cycles with the added conditions that each of $1,2,\dots,l-1$ appears in the first cycle and each of $l,l+1,\dots,k$ appears in the second cycle.  (As usual, we assume that the cycles are arranged with their smallest elements first.)  Indeed, if the first cycle has length $t-1$, then we have
$${t-2 \choose l-2} \cdot (l-2)!$$ ways to position $2,3,\dots,l-1$; $${n-k \choose t-l} \cdot (t-l)!$$ ways to choose and position the remaining elements in the first cycle; $${n-t \choose k-l}\cdot (k-l)!$$ ways to position $l+1,l+2,\dots,k$ in the second cycle; and $$(n-t-k+l)!$$ ways to arrange the remaining elements in the second cycle.  As $t$ ranges from $l$ to $n-k+l$, we account for all desired cycle-decompositions, verifying the left-hand side of (\ref{Stirlingl>=2'}).

To arrive at the right-hand side, we start by arranging the elements $2,3,\dots,l-1$ in the first cycle (all behind 1) and arranging $l+1,l+2,\dots,k$ in the second cycle (all behind $l$); there are
$$(l-2)! \cdot (k-l)!$$ ways to do this.  The elements $k+1,k+2,\dots,n$ can then be inserted, one at a time, to the right of previously-placed elements; this can be accomplished in $$\frac{(n-1)!}{(k-1)!}$$ ways.  The result now follows.

Finally, we turn to the case of $l=1$ and $d=2$; from (\ref{Vformula}) we have
\begin{eqnarray*} 
V_{n,2}(k,1) & = & \sum_{t=2}^{n-k+1} \sum_{s=1}^{\min\{t-1,2\}} P^{\surd}_n(k,1;t,s) \cdot (n+1-s) + 
    \sum_{s=2}^{2} P^{X}_n(k,1;s) \cdot (n+1-s) \\ \\
& = & \sum_{t=2}^{n-k+1} \sum_{s=1}^{2} P^{\surd}_n(k,1;t,s) \cdot (n+1-s) - P^{\surd}_n(k,1;2,2) \cdot (n-1) + 
    P^{X}_n(k,1;2) \cdot (n-1),
\end{eqnarray*}
which, using (\ref{Psurd}) and (\ref{PX}) and simplifying, becomes
\begin{eqnarray*}
V_{n,2}(k,1)=\frac{2n-1}{{n \choose k}} \cdot \sum_{t=2}^{n-k+1} {n-t \choose k-1} \cdot \frac{1}{t-1} - \frac{k}{n} \cdot (n-k-1).  
\end{eqnarray*}
Our claim for $V_{n,2}(k,1)$ now follows from (\ref{Stirlingl=1}).
$\Box$

\begin{thm}
For a fixed value of $d$, let $$c_d= \lim_{n \to \infty} \frac{\max\{V_{n,d}(k,l) \mid k, l \}}{n}.$$  Then $c_1=1/\mathrm{e} \approx 0.36788,$  and $c_2=2x-x^2 \approx 0.51239$ where $x$ is the smaller root of the equation $$2x - 2\log x = 3.$$

\end{thm}

{\em Proof.}  We will use Theorem \ref{Vnd} to find the values of $k$ and $l$ for which $V_{n,d}(k,l)$ is maximal for $d=1$ and $d=2$.

Consider first $d=1$, in which case for $l=1$ we have $$V_{n,1}(k,1)=k \cdot (H_n-H_k)=k \cdot (n-k) \cdot \mathrm{Average}\left\{\frac{1}{n-i} \; | \; i=1,2,\dots,n-k\right\}$$ and for $l \geq 2$ we have
\begin{eqnarray*}
V_{n,1}(k,l) & = &  \frac{n}{(l-1) \cdot {n \choose k}} \cdot \left({n-1 \choose k-1} - {n-l \choose k-l} \right) \\
& = & \frac{n}{(l-1) \cdot {n \choose k}} \cdot \left({n-2 \choose k-1} + {n-3 \choose k-2 } + \cdots + {n-l \choose k-l+1} \right) \\
& = & \frac{n}{(l-1) \cdot {n \choose k}} \cdot \left({n-2 \choose n-k-1} + {n-3 \choose n-k-1 } + \cdots + {n-l \choose n-k-1} \right) \\
& = & \frac{n}{ {n \choose k}} \cdot \mathrm{Average}\left\{ {n-j \choose n-k-1}  \; | \; j=2,3,\dots,l \right\}
\end{eqnarray*}
where $\mathrm{Average} \; S$ is the arithmetic average of the finite set of real numbers $S$.  Note that $\frac{1}{n-i}$ is an increasing function of $i$ but ${n-j \choose n-k-1}$ is a decreasing function of $j$, thus we have
$$V_{n,1}(k,1) \geq k \cdot (n-k) \cdot \frac{1}{n-1}$$ and
$$V_{n,1}(k,l) \leq \frac{n}{ {n \choose k}} \cdot {n-2 \choose n-k-1} = k \cdot (n-k) \cdot \frac{1}{n-1} $$
as $l \geq 2$.  Therefore, for fixed $n$ and $k$, $V_{n,1}(k,l)$ is maximal for $l=1$.

To find the maximum of $V_{n,1}(k,1)$, we let $x$ denote the limit of $k/n$ as $n$ approaches infinity; we then have  
$$V_{n,1}(k,1)= k \cdot (H_n-H_k) \sim  - n \cdot x \cdot \log x,$$ and this attains its maximum when $x =1/ \mathrm{e}$.  Therefore, for $d=1$, $V_{n,1}(k,l)$ attains its maximum when $l=1$ and $k \sim n/\mathrm{e}$, in which case we have
$V_{n,1}(k,1) \sim  n/ \mathrm{e}.$

The situation is somewhat more complicated if $d=2$.  We still find that for $l \geq 2$, 
$$V_{n,2}(k,l)=\frac{2n-1}{ {n \choose k}} \cdot \mathrm{Average} \left\{ {n-j \choose n-k-1}  \; | \; j=2,3,\dots,l \right\}$$ which is a decreasing function of $l$, but comparing $V_{n,2}(k,1)$ and $V_{n,2}(k,2)$ is not possible.  We will, however, show that, as $n$ approaches infinity, the maximum value of $V_{n,2}(k,1)$ is more than the maximum value of $V_{n,2}(k,2)$, as follows.

Clearly, $$V_{n,2}(k,2)=\frac{2n-1}{ {n \choose k}} \cdot {n-2 \choose n-k-1}=\frac{2n-1}{n^2-n} \cdot k \cdot (n-k)$$ is maximal when $k=n/2$; in which case we have $V_{n,2}(k,2) \sim n/2$.

To find the maximum of $V_{n,2}(k,1)$, we let $x$ denote the limit of $k/n$ as $n$ approaches infinity, as before; we then have  
$$V_{n,2}(k,1)= \frac{2n-1}{n} \cdot k \cdot (H_n-H_k) - \frac{k}{n} \cdot n \cdot \left(1-\frac{k}{n} - \frac{1}{n} \right) \sim  - 2n \cdot x \cdot \log x - n \cdot x \cdot (1-x).$$ This attains its maximum when $2x - 2\log x = 3$ or $x \approx 0.30171$, in which case $$V_{n,2}(k,1) \sim n \cdot x \cdot (x-2\log x -1)=n \cdot x \cdot (2-x)$$ or about $0.51239 n $.  Since this is larger than $n/2$, we have established our claim.  $\Box$

\end{document}